\documentclass[11pt]{amsart}
\usepackage[tcidvi]{graphicx}
\usepackage{amscd}
\usepackage{amsmath}
\usepackage{amsfonts}
\usepackage{amssymb}
\theoremstyle{plain}

\numberwithin{equation}{section}

\newtheorem{thmintro}{Theorem}
\newtheorem{lemintro}[thmintro]{Lemma}
\newtheorem{defintro}[thmintro]{Definition}
\newtheorem{theorem}{Theorem}[section]
\newtheorem{corollary}[theorem]{Corollary}
\newtheorem{lemma}[theorem]{Lemma}

\theoremstyle{definition}
\newtheorem{definition}[theorem]{Definition}

\newtheorem*{remark}{Remark}%
\begin{document}
\large

\title[The Continuous Wavelet Transform]{The Continuous Wavelet Transform and Symmetric Spaces}
\author{R. Fabec and G. \'Olafsson}
\address{Department of Mathematics, Louisiana State University, Baton Rouge,
LA 70803, USA}
\email{fabec@math.lsu.edu}
\address{Department of Mathematics, Louisiana State University, Baton Rouge,
LA 70803, USA}

\email{olafsson@math.lsu.edu}
\thanks{Research by G. \'Olafsson supported
by NSF grant DMS-0070607 and the MSRI}
\begin{abstract}
The continuous wavelet transform has become a widely used tool in 
applied science during the last decade. In this article we discuss 
some generalizations coming from actions of closed subgroups of 
$\mathrm{GL}(n,\mathbb{R})$ acting on $\mathbb{R}^n$. In particular, 
we propose a way to invert the wavelet transform in the case where the 
stabilizer of a generic point in $\mathbb{R}^n$ is not compact,
but a symmetric subgroup, a case 
that has not previously been discussed in the literature.  
\end{abstract}

\maketitle

\section*{Introduction, Wavelets on $ax+b$-group}

The continuous wavelet transform has become a widely used tool in 
applied science during the last decade. The best known example deals 
with wavelets on the real line. Each wavelet by taking ``matrix 
coefficients'', i.e., inner products with translations and dilations 
of the wavelet, defines a wavelet transform which can be used to 
reconstruct the function from the dilations and translations of the 
wavelet. This process helps compensate for the local--nonlocal 
behavior of the Fourier transform.

Translations and dilations form the so-called $ax+b$-group of 
transformations of the real line.  These act in a natural manner
as unitary operators on $L^2(\mathbb{R})$.   
In that way, wavelet transforms are simply a part of the 
representation theory of the $ax+b$-group.  This observation is 
the basis for the generalization of the continuous wavelet transform 
to higher dimensions and more general settings.  There have 
been further attempts to generalize these ideas to arbitrary groups; 
see \cite{FM01} and the references therein. 
In this article we review some of the basic ideas 
for the general wavelet transform for groups acting on $\mathbb{R}^n$. 
We start by reviewing the classical one dimensional wavelet transform. 
In this simple setting, the usual definition of a wavelet is 
equivalent to the corresponding matrix coefficients being square 
integrable on the $ax+b$--group. A natural generalization of this to an 
arbitrary representation $(\pi ,\mathbb{H})$ of a topological group 
$G$ is to say, that a vector $u\in\mathbb{H}\setminus \{0\}$ is a 
wavelet if $a\mapsto (v\mid\pi (a)u)$ is square integrable for all 
$v\in\mathbb{ H}$. The generalized wavelet transform is then 
$W_v(a):=(v\mid\pi (a)v)$. The basic facts for this transform and, in 
particular, the inversion formula are presented in section \ref{sec1}.

In section \ref{sec2} these notions are applied to topological groups $H$ 
acting on $\mathbb{R}^n$ by a representation $\mathbf{\pi}$. Let 
$G=H\times_{\mathbf{\pi}}\mathbb{R}^n$ be the semi-direct product of 
$H$ and $\mathbb{R}^n$. Then $G$ acts on $\mathbb{R}^n$ and $L^2 
(\mathbb{R}^n)$. When the action of $H$ on $\mathbb{R}^n$ is 
sufficiently well behaved, the decomposition of $L^2(\mathbb{R} ^n)$ 
under $G$ can be described in terms of the orbits in $\mathbb {R}^n$ 
under the transpose (contragredient) representation $\pi'$ where 
$\pi'(a)=\pi (a^{-1})^t$. In general, one may not even have 
irreducible subrepresentations inside $L^2(\mathbb R^n)$, but in the 
case when there are open orbits, one has sufficiently many irreducible 
subspaces to decompose all function in $L^2(\mathbb R^n)$. In this 
situation, those irreducible subspaces corresponding to orbits which 
have compact stabilizers will be precisely those subspaces for which 
wavelets exist. Some typical examples come from identifying the 
nilradical of a parabolic subalgebra $\frak{p}=\frak{m}\oplus \frak
{a}\oplus \frak{n}$ in a semisimple Lie algebra $\frak{g}$ with 
$\mathbb{R}^n$. Several such examples can be found in the literature
\cite{BR}.

Thus an interesting situation occurs when $\mathbb{R}^{n}$ decomposes 
up to set of measure zero into finitely many open orbits $\left\{
\mathcal{O}_{i}\right\}_{i=1}^{k}$ under $\mathbf{\pi}(H)^{t}$. In 
this case, $L^{2}(\mathbb{R}^{n})=\bigoplus_{i=1}^{k}L_{\mathcal{O}_ 
{i}} ^{2}(\mathbb{R}^{n})$ where $L_{\mathcal{O}_{i}}^{2}(\mathbb{R}^ 
{n})=\left\{ f\in L^{2}(\mathbb{R}^{n})\mid\hat{f}|_{\mathbb{R}^{n}
\setminus\mathcal{O} _{i}}=0\right\}$ gives a decomposition of $L^
{2}(\mathbb{R}^{n})$ into irreducible subspaces. As mentioned earlier 
and as will be seen in section \ref{sec2}, the subspace $L_{\mathcal 
{O}_{i}}^ {2}(\mathbb{R}^{n})$ contains a wavelet vector in the above 
sense if and only if the stabilizers of points in $\mathcal{O}_{i}$ 
are compact. In each of these cases, there will be wavelet transforms. 
The question then becomes what can be done in the general case where 
the stabilizer of a point in the oribit $\mathcal{O}_i$ is 
non-compact. This situation occurs already in the simple example of 
the group $\mathbb{R}^+\mathrm{SO}_o(1,n)$ acting on $\mathbb{R}^n$ 
which we discuss in detail in section \ref{sec4}. It falls into the 
category where one has symmetric orbits, a case which we can answer.

Let $x\in\mathcal{O}_i$ and $L=\left\{a\in H\mid\mathbf{\pi}(a) 
^tx=x\right\}$. Assume that the group $H$ is reductive and that there 
exists an involution $\tau :H\rightarrow H$ such that $H_o^ {\tau}
\subset L\subset H^ {\tau}=\left\{a\in H\mid\tau (a) =a\right\}$. Then 
using a classical result of Matsuki \cite{Ma79,GO87}, one can show 
that there exist a closed group $R=R_i\subset H$ and a finite set 
$x_1\ldots ,x_r\in H$ such that each $ Rx_q$ is open, $Rx_q\cap Rx_
{q'}=\emptyset$ ($q\not =q'$), $Rx_1\cup\ldots\cup Rx_r$ is dense in 
$\mathcal{O}_i$, and the stabilizer of each $x_q$ in $R$ is compact 
(cf. \cite{GO87}). These ideas are used in section \ref{sec3} to 
decompose each of the $G$ irreducible spaces $L_{\mathcal{O}_i}^2
(\mathbb{R}^n)$ into irreducible subspaces under $G_i=R_i\times_{\pi}
\mathbb R^n$ having $G_i$-wavelet vectors. These vectors and the 
results from section \ref{sec2} can then be used to construct a 
wavelet transform and an inversion formula for the wavelet transform. 
We give simple examples of this in the final section.

The $ax+b$ group is the group of affine transformations $T_{a,b}$ 
acting on $\mathbb{R}$ by dilation by $a>0$ followed by translation by 
$b$. Thus
\[
T_{a,b}x=ax+b.
\]
Note
\[
T_{a,b}T_{a^{\prime},b^{\prime}}x=aa^{\prime}x+ab^{\prime}+b=T_{aa^{\prime
},b+ab^{\prime}}
\]
and
\[
T_{a,b}^{-1}=T_{\frac{1}{a},-\frac{b}{a}}.
\]

We identify $(a,b)$ where $a>0$ and $b\in\mathbb{R}$ to the transformation
$T_{a,b}$ and take $G$ to be collection of all such pairs. The action of $G$
on $\mathbb{R}$ ``induces'' a natural action $\mathbf{\rho}$ of $G$ on the
space of $L^{2}$ functions on $\mathbb{R}$. This action is called the left
\textit{quasi-regular representation} of $G$ on $\mathbb{R}$. It is defined
by
\[
\rho(a,b)f(x)=J(T_{a,b}^{-1})^{\frac{1}{2}}f(T_{a,b}^{-1}x)=a^{-\frac{1}{2}
}f(\frac{x-b}{a})
\]
where $J(\cdot)$ denotes the Jacobian of the transformation. If we 
allow $a$ also to be negative the action has to be replaced by $\left| 
a\right| ^{-1/2}f\left( \frac{x-b}{a}\right)$. The presence of $J$ is 
to force the natural action to be unitary. Because of the composition 
laws in $G$, $\rho(a,b)\rho(a^{\prime},b^{\prime})=\rho((a,b)(a^
{\prime},b^{\prime}))$. The space $L^{2}(\mathbb{R})$ under the action 
$\rho$ is known to contain precisely two closed proper invariant 
subspaces. Thus to decompose a function $f\in L^{2}(\mathbb{R})$ into 
``wavelets'', one needs to know how one can project $f$ into each of 
these subspaces. The wavelet transform and the square integrability of 
$\rho$ provides the key ingredient. For completeness, we describe the 
decomposition of $L^{2}(\mathbb{R})$ into its two closed proper 
invariant subspaces. The argument we give is elementary; it decomposes 
$L^{2}(\mathbb{R})$ into a the sum $H_{+}^2\oplus H_{-}^2$, the 
classical Hardy spaces of functions $f$ satisfying $\hat{f}(\pm y)=0$ 
for $y<0$.

\begin{thmintro} 
The two subspaces $H_{+}^{2}$ and $H_{-}^{2}$ are the only proper 
closed invariant subspaces under $\rho$ and $L^{2}(\mathbb{R})=H_{+}^
{2}\oplus H_{-}^{2}$.
\end{thmintro}

\begin{proof}
The Fourier transform $f\mapsto\hat{f}$ given formally by
\[
\hat{f}(y)=\frac{1}{\sqrt{2\pi}}\int f(x)e^{-ixy}\,dx
\]
is a linear isometry from $L^{2}(\mathbb{R})$ onto $L^{2}(\mathbb{R})
$. If one defines $\hat{\rho}(a,b)$ by
\[
\hat{\rho}(a,b)\hat{f}=\widehat{\,\rho(a,b)f\,},
\]
then
\begin{equation}\tag{1}
\begin{split}
\hat{\rho}(a,b)\hat{f}(y) &  =\frac{1}{\sqrt{2\pi}}\int a^{-1/2}
f(a^{-1}(x-b))e^{-ixy}\,dx\\
&  =\frac{1}{\sqrt{2\pi}}\int a^{-\frac{1}{2}}f(a^{-1}x)e^{-i(x+b)y}\,dx\\
&  =\frac{1}{\sqrt{2\pi}}e^{-iby}\int a^{-\frac{1}{2}}f(x)e^{-iaxy}\,d(ax)\\
&  =\frac{1}{\sqrt{2\pi}}a^{\frac{1}{2}}e^{-iby}\int f(x)e^{-iaxy}dx\\
&  =\sqrt{a}e^{-iby}\hat{f}(ay).
\end{split}
\label{rhohat}
\end{equation}
Since $a>0$, this implies the subspaces $\hat{H}_{+}^{2}$ and $\hat{H}_{-}
^{2}$ are invariant under $\hat{\rho}$ and $\hat{H}_{+}^{2}\oplus \hat{H}_{-}
^{2}=\hat{L}_{2}(\mathbb{R})$.

Let $S$ be a closed invariant subspace of $\hat{L}_{2}(\mathbb{R})$. Note if
$h\in L^{1}(\mathbb{R})$ and $\hat{f}\in S$, then the vector
\[
 M_{\hat{h}}\hat{f}(y)\equiv\frac{1}{\sqrt{2\pi}}\int h(b)\rho(1,
 b)\hat {f}(y)\,db=\frac{1}{\sqrt{2\pi}}\int h(b)e^{-iby}\hat{f}
 (y)\,db=\hat{h}(y)\hat {f}(y)
\]
is a limit in $\hat{L}_{2}(\mathbb{R})$ of linear combinations of 
vectors of form $\rho(1,b_{k})f$ and thus is in $S$. Hence $\hat{h}
\hat{f}\in S$ for any $h\in L^{1}(\mathbb{R})$. Since $|\hat{h}|_
{\infty}\leq|h|_{1}$, one can argue that given any $\hat{f}\in S$ and 
any continuous bounded function $h$, one has $h\hat{f}\in S$. 
Consequently $h\hat{f}\in S$ for any bounded measurable function $h$. 
In particular, if $P$ is an orthogonal projection commuting with all 
the operators $\hat{\rho}(1,b)$, $P$ commutes with all multiplication 
operators $M_{h}$ and thus must have form $M_{\chi_{E}}$ where $E$ is 
some measurable set; i.e., $Pf=\chi_{E}f$ for all $f\in\hat{L}_{2}
(\mathbb{R})$. If $P$ is the orthogonal projection onto $S$, we see 
since
\begin{align*}
\sqrt{a}\chi_{E}(ay)f(ay) &  =\sqrt{a}Pf(ay)\\
&  =\hat{\rho}(a,1)P\hat{f}(y)\\
&  =P\hat{\rho}(a,1)\hat{f}(y)\\
&  =\chi_{E}(y)\sqrt{a}\hat{f}(ay),
\end{align*}
that $\chi_{a^{-1}E}=\chi_{E}$ for all $a>0$. Thus $E$ is essentially 
the empty set, the set $(-\infty,0]$, the set $[0,\infty)$, or 
$\mathbb{R}$. This shows $L^{2}(\mathbb{R})$ has precisely four closed 
invariant subspaces, the only proper ones being $H_{+}^{2}$ and $H_
{-}^{2}$.
\end{proof}

The following discussion on wavelets can be extended to all of $L^{2} 
(\mathbb{R})$, but we will concentrate on the Hardy space $H_{+}^{2}$, 
for the analysis is the same in the other case. We normalize the Haar 
measure on $G$ to be
\[
dg=\frac{dadb}{2\pi a^{2}}.
\]

\begin{defintro}
Let $\psi$ be a nonzero vector in $H_{+}^{2}$. Then $\psi$ is called a
\textbf{wavelet} if $\hat{\psi}\in L^{2}(\mathbb{R},\frac{dt}{|t|}).$
\end{defintro}

Let $\psi$ be a wavelet. By (\ref{rhohat}), $\rho(a,b)\psi$ is a 
wavelet for each $(a,b)$. If $f\in L^{2}(\mathbb{R})$, then the 
\textit{wavelet transform} of $f$ is defined by
\[
W_{\psi}f(a,b):=(f\mid\rho(a,b)\psi)=a^{-1/2}\int_{-\infty}^{\infty
}f(t)\overline{\psi\left(  \frac{t-b}{a}\right)  }\,dt.
\]
The function $G\ni (a,b)\mapsto Wf(a,b)$ is  continuous on $G$.

Note if $\psi$ is any nonzero element in $H_{+}^{2}$, then the 
vectors $\rho(a,b)\psi$ span a dense subspace of $H_{+}^{2}$. It 
is well known that any function $f$ can be recovered from the 
wavelet transform $W_{\psi}f$. The key is the square integrability 
of the wavelet transform.

\begin{lemintro}
Let the notation be as above. Then
\[W_{\psi}f(a,b)=a^{1/2}\int\hat{f}(\omega)\overline{\hat{\psi}
(a\omega )}\,e^{i\omega b}\,d\omega\,.
\]
\end{lemintro}

\begin{proof}
Note by (\ref{rhohat}),
\begin{align*}
W_{\psi}f(a,b) &  =(f,\rho(a,b)\psi)\\
&  =(\hat{f},\hat{\rho}(a,b)\hat{\psi})\\
&  =\int\hat{f}(y)a^{1/2}\overline{e^{-iby}\hat{\psi}(ay)}\,dy\\
&  =a^{1/2}\int\hat{f}(y)\overline{\hat{\psi}(ay)}\,e^{iby}\,dy
\end{align*}
\end{proof}

\begin{thmintro}
\label{t63} Let $\psi$ be a wavelet. Set $C_{\psi}^{2}=\int_{0}^{\infty}
\frac{|\hat{\psi}(\omega)|^{2}}{\omega}\,d\omega>0$. Then for any $f\in
H_{+}^{2}$, $W_{\psi}f$ is in $L^{2}(G)$. Moreover, for $f,h\in H_{+}^{2}$,
\[(W_{\psi}f,W_{\psi}g)=C_{\psi}^{2}(f,g)\,.
\]
\end{thmintro}

Recall that if $(\pi,\mathbb{H})$ and $(\rho,\mathbb{K})$ are two 
representations of a group $G$, then a continuous map $T:\mathbb 
{H}\rightarrow\mathbb{K}$ is called an \textit{intertwining 
operator} if for all $x\in G$ we have $T\pi(x)=\rho(x)T$. We can 
then describe the wavelet transform $W_{\psi}$ as an intertwining 
operator $H_{+}^{2}\rightarrow L^{2}(G)$:

\begin{lemintro}
Let $\psi\in H_{+}^{2}$ be a wavelet. Then
\[
    H_{+}^{2}\ni f\mapsto W_{\psi}f\in L^{2}(G)
\]
is continuous with norm $C_{\psi}>0$. Furthermore, $W_{\psi}$ 
intertwines the representation $\rho$ and the left regular 
representation $L$ on $G$; i.e., $W_{\psi}\rho(a,b)=L(a,b)W_
{\psi}$ for all $(a,b)\in G$.
\end{lemintro}

This lemma is the basis for understanding and generalizing wavelet 
transforms to other groups. The reproducing property and inversion 
formula for the wavelet transform are simple consequences of 
intertwining the representation with the regular representation and 
the irreducibility of the representation on $H_{+}^2$. We discuss this 
general framework in the next section.  

\section{Wavelets and Square Integrable Representations}\label{sec1}

In this section $G$ denotes a locally compact Hausdorff topological 
group. We fix a \textit{left} Haar measure $d\mu=dx$ for $G$.
We define a \textit{wavelet vector} with respect to
an unitary representation of $G$ in the following way:

\begin{definition}
Suppose $(\rho,\mathbb{H})$ is an irreducible unitary representation 
of $G$. Then a \textbf{wavelet}
$\psi$ for $\rho$ is a nonzero vector such that 
$\rho_{v,\psi }(g):=(v\mid\rho(g)\psi)$ is in $L^{2}(G)$ for all 
$v\in\mathbb{H}$.
\end{definition}

\begin{lemma}
Assume $\psi$ is a wavelet for the irreducible representation $\rho$ 
on the Hilbert space $V_\rho$.  Then the map $W_{\psi}:V_{\rho}
\rightarrow L^{2}(G)$ sending $v$ to $\rho_{v,\psi}$ is a continuous 
intertwining operator.
\end{lemma}

\begin{proof}
This linear mapping has closed graph. Hence it is continuous. We  
show that $W_{\psi}$ is intertwining. Indeed,
\[
W_{\psi}(\rho(x)u)(y)=(\rho(x)u\mid\rho(y)\psi)=(u\mid\rho(x^{-1}y)
\psi )=L_{x}W_{\psi}u(y)\,.
\]
\end{proof}

\begin{remark}
Duflo and Moore \cite{DM76} showed for any $\psi$, $\rho_{u,\psi}$ is 
square integrable for all $u$ if it is square integrable for one 
nonzero $u$ and defined an irreducible representation to be square 
integrable if it has a nonzero square integrable matrix coefficient. 
Hence the wavelets for $\rho$ are the $\psi$ which give square 
integrable matrix coefficients. Moreover, since $||\rho_{u,\rho(g)
\psi}||_{2}^{2}=\Delta(g)^{-1}\,||\rho_{u,\psi} ||_{2}^{2}$, the 
wavelets form an invariant and hence dense linear subspace of 
$V_{\rho}$. 

The special case where $G=P=HN$ is a \textit{parabolic subgroup} in a 
semisimple Lie group is discussed in detail in \cite{LW78,LW82,W79}. 
Here $H$ is the Levi-factor of $P$ and $N$ is the nilradical. In 
particular, $H$ is reductive and acts on the Lie algebra of $N$, which 
in many situations is isomorphic to $\mathbb{R}^n$. These special 
cases are important for our discussions on the generalization of the 
wavelet transform.

\end{remark}

\begin{definition}
Let $X$ be a set and let $\mathbb{H}$ be a Hilbert space of functions 
$f:X\rightarrow\mathbb{V}$ where $\mathbb{V}$ is a Hilbert space. The 
space $\mathbb{H}$ is called a \textbf{reproducing Hilbert space} if 
the evaluation maps $\mathrm{ev}_{x}:\mathbb{H}\rightarrow\mathbb{V}$, 
$f\mapsto f(x)$, are continuous for all $x\in X$.
\end{definition}

\begin{lemma}
Let $\mathbb{H}$ be a reproducing Hilbert space. Define $K(x,y):
=\mathrm{ev} _{x}\mathrm{ev}_{y}^{\ast}:\mathbb{V}\rightarrow\mathbb
{V}$. Let $K_{y}(x):=K(x,y)$. Then the following hold:
\begin{enumerate}
\item\vskip -.2cm 
$K(x,y):\mathbb{V}\rightarrow\mathbb{V}$ is continuous and linear; 
i.e., $K(x,y)\in\mathrm{Hom}(\mathbb{V},\mathbb{V})$

\item $K_{y}(\cdot)v\in\mathbb{H}$ for all $v\in\mathbb{V}$, and
$(f\mid K_{y}v)=(f(y)\mid v)$ for all \ $f\in\mathbb{H}$ and all 
$v\in\mathbb{V}$.

\item 
The set of 
finite linear combinations $\sum c_{j}K_{y_{j}}(\cdot)v_{j}$ is dense 
in $\mathbb{H}$.

\item 
$K(x,y)^{\ast}=K(y,x)$
\end{enumerate}
\end{lemma}

\begin{proof}
As $\mathrm{ev}_{x}:\mathbb{H}\rightarrow\mathbb{V}$ is continuous and 
linear for all $x\in X$, it follows that $\mathrm{ev}_{x}^{\ast}:
\mathbb{V}
\rightarrow\mathbb{H}$ is continuous and linear. Hence $K(x,y):\mathbb{V}
\rightarrow\mathbb{V}$ is continuous and linear. We have $K(\cdot
,y)v=\mathrm{ev}_{y}^{\ast}v\in\mathbb{H}$ for all $y\in\mathbb{V}$.

Let $u\in\mathbb{V}$ and $y\in X$. Then
\begin{align*}
(f(y)   \mid v) & =(\mathrm{ev}_{y}f\mid v)\\
&  =(f\mid\mathrm{ev}_{y}^{\ast}v)\\
&  =(f\mid K(\cdot,y)v)\,.
\end{align*}
Hence (b) follows. Assume that $f\in\mathbb{H}$ is perpendicular to 
all the linear combinations $\sum c_{j}K_{y_{j}}(\cdot)v_{j}$. Then in 
particular
\[
(f(y)\mid v)=(f\mid K(\cdot,y)v)=0
\]
for all $y\in X$ and all $v\in\mathbb{V}$. Hence $f(y)=0$ for each $y$ 
and thus $f=0$. So (c) follows. Finally (d) follows directly from the 
definition as
\[
K(x,y)^{\ast}=\left( \mathrm{ev}_{x}\mathrm{ev}_{y}^{\ast}\right) ^
{\ast }=\mathrm{ev}_{y}\mathrm{ev}_{x}^{\ast}=K(y,x)\,.
\]
\end{proof}

\begin{definition}
The map $K:X\times X\rightarrow\mathrm{Hom}(V,V)$ is called the
\textbf{reproducing kernel} of $\mathbb{H}$ .
\end{definition}

Let $(\pi,\mathbb{H})$ be a unitary representation of $G$ and let $f$ 
be a measurable function on $G$ such that the integral $\int_{G}f
(x) (\pi(x)u\mid v)\,dx$ converges for all $u,v\in\mathbb{H}$ and for 
which there is a constant $C$ with $\left| \int_{G} f(x)(\pi(x)u\mid 
v)\,dx\right| \leq C||u||\,||v||$. Then we can define a continuous 
linear map $\pi(f):\mathbb{H}\rightarrow\mathbb{H}$ by
\[
(\pi(f)u\mid v)=\int_{G}f(x)(\pi(x)u\mid v)\,dx\,.
\]

\begin{theorem}\label{Theorem 1.6}
Let $(\rho,\mathbb{H})$ and $(\tau,\mathbb{K})$ be square integrable 
representations with wavelet vectors $\psi$ and $\eta$, respectively. 
Then the following hold:
\begin{enumerate}
\item\vskip -.2cm  
If $\rho$ and $\tau$ are not equivalent, then $(W_{\psi}u\mid W_{\eta 
}v)=0$ for all $u\in V_{\rho}$ and all $v\in V_{\tau}$. Thus the 
images of $W_{\psi}$ and $W_{\eta}$ are orthogonal subspaces of $L^
{2}(G)$.

\item  
Let $f\in L^{2}(G)$. Then $\rho(f)\psi$ is (weakly) defined and 
$W_{\psi}^{\ast}f=\rho(f)\psi$.

\item  
If $\rho=\tau$, then there exists a constant $C_{\psi,\eta}$ such that
\[
(W_{\psi}u\mid W_{\eta}v)=C_{\psi,\eta}(u\mid v)\,
\]
for all $u,v\in V_{\rho}$.

\item 
$\mathrm{Im}(W_{\psi})$ is a closed reproducing Hilbert subspace of 
$L^{2}(G)$ contained in $L^{2}(G)\cap C(G)$.

\item  We have $C_{\psi,\psi}>0$.
Let $c_{\psi}=1/\sqrt{C_{\psi,\psi}}>0$. Define $U_{\psi}:\mathbb{H}
\rightarrow L^{2}(G)$ by $u\mapsto c_{\psi}W_{\psi}u=W_{c_{\psi}\psi}u$. 
Then $U_{\psi}$ is an unitary isomorphism $\mathbb{H}\simeq\mathrm
{Im}(W_{\psi})$.

\item  
Assume that $\rho=\tau$ and that $C_{\psi,\eta}\not =0$. Then
\[
u=\frac{1}{C_{\psi,\eta}}\rho(W_{\psi}u)\eta=\frac{1}{C_{\psi,\eta}}
\int _{G}(u\mid\rho(x)\psi)\rho(x)\eta\,dx,
\]
(weakly) for all $u\in V_{\rho}$. In particular,
\[
u=\frac{1}{C_{\psi,\psi}}\rho(W_{\psi}u)\psi=\frac{1}{C_{\psi,\psi}}\int
_{G}(u\mid\rho(x)\psi)\rho(x)\psi\,dx\,.
\]

\item  
The reproducing kernel for $\mathrm{Im}(W_{\psi})$ is given by $K(x,
y)=C_{\psi,\psi}^{-1}W_{\psi}\psi(y^{-1}x)$

\item  
If $f\in\mathrm{Im}(W_{\psi})$, then $f=C_{\psi,\psi}^{-1}f\ast W_
{\psi}(\psi)$.
\end{enumerate}
\end{theorem}

\begin{proof}
Define $\beta:\mathbb{H}\times\mathbb{K}\rightarrow\mathbb{C}$ by
\[
\beta(u,v):=(W_{\psi}u\mid W_{\eta}v)=\int W_{\psi}u(x)\,\overline{W_{\eta
}v(x)}\,dx\,.
\]
As the wavelet transform is continuous, it follows that
\[
|\beta(u,v)|=|(W_{\psi}u\mid W_{\eta}v)|\leq||W_{\psi}u||\,||W_{\eta}
v||\leq||W_{\psi}||\,||W_{\eta}||\,||u||\,||v||\,.
\]
Consequently, there is a continuous linear map $T:\mathbb{H}\rightarrow
\mathbb{K}$ such that
\[
\beta(u,v)=(Tu\mid v)\,.
\]
Let $a\in G$. Then, using that the wavelet transform is intertwining,
\begin{align*}
(T\rho(a)u\mid v)  &  =\beta(\rho(a)u,v)\\
&  =\int[W_{\psi}\rho(a)u](x)\overline{W_{\eta}v(x)}\,dx\\
&  =\int W_{\psi}u(a^{-1}x)\overline{W_{\eta}v(x)}\,dx\\
&  =\int W_{\psi}u(x)\overline{W_{\eta}v(ax)}\,dx\\
&  =\int W_{\psi}u(x)\overline{[W_{\eta}\tau(a^{-1})v](x)}\,dx\\
&  =\beta(u,\tau(a^{-1})v)\\
&  =(Tu\mid\tau(a)^{\ast}v)\\
&  =(\tau(a)Tu\mid v)\,.
\end{align*}

As this holds for all $u$ and $v$, it follows that $T\rho(a)=\tau(a)
T$. Thus $T$ is an intertwining operator. It follows by Schur's Lemma 
that $T$ is zero when $\rho$ is not equivalent to $\tau$, and 
$T=C_{\psi,\eta}I$ for some scalar when $\rho=\tau$. Moreover, if 
$\psi=\eta$, then $T>0$ and thus $C_{\psi,\psi}>0$. Consequently, 
$\frac{1}{\sqrt{C_{\psi,\psi}}}W_{\psi}$ is an isometry of $V_{\rho}$ 
onto its image in $L^{2}(G)$. Hence (a), (c), and (e) hold.

Let $f\in L^{2}(G)$ and $u\in V_{\rho}$. Then
\begin{align*}
(W_{\psi}^{\ast}f\mid u)  &  =(f\mid W_{\psi}u)\\
&  =\int f(x)\,\overline{(u\mid\rho(x)\psi)}\,dx\\
&  =\int f(x)(\rho(x)\psi\mid u)\,dx\\
&  =(\rho(f)\psi\mid u)
\end{align*}
Hence $\rho(f)\psi$ is weakly defined and (b) holds.

(d) The functions $f=W_{\psi}u\in\operatorname{Im}(W_{\psi})$ are continuous
functions. In particular, $\mbox{\rm ev}_{a}(f)=f(a)$ is well defined on
$\mbox{\rm Im}(W_{\psi})$. Since $U_{\psi}$ is an isometry,
\begin{align*}
|\mbox{\rm ev}_{a}(f)|  &  =|W_{\psi}u(a)|\\
&  =|(u,\rho(a)\psi)|\\
&  \leq||u||\,||\psi||\\
&  \leq||U_{\psi}u||\,||\psi||\\
&  =c_{\psi}||W_{\psi}u||\,||\psi||\\
&  =c_{\psi}||\psi||\,||f||_{2}.
\end{align*}
Hence point evaluations are continuous linear mappings.

Note (f) follows from (b) and (c) for:
\begin{align*}
\frac{1}{C_{\psi,\eta}}(\rho(W_{\psi}u)\eta\mid v)&=
\frac{1}{C_{\psi,\eta}}(W_{\eta}^{\ast}(W_{\psi}u)\eta\mid v)\\
&  =\frac{1}{C_{\psi,\eta}}(W_{\psi}u\mid W_{\eta}v)\\
&  =(u\mid v).
\end{align*}

Note if $f=W_{\psi}u$, then
\begin{align}
f(y)  &  =W_{\psi}u(y)\nonumber\\
&  =(u\mid\rho(y)\psi)\nonumber\\
&  =\frac{1}{C_{\psi,\psi}}(\rho(W_{\psi}u)\psi\mid\rho(y)\psi)\nonumber\\
&  =\frac{1}{C_{\psi,\psi}}\int W_{\psi}u(x)(\rho(x)\psi\mid\rho
(y)\psi)\,dx\nonumber\\
&  =\frac{1}{C_{\psi,\psi}}\int f(x)(\rho(y^{-1}x)\psi\mid\psi
)\,dx\label{eq1.2}\\
&  =\frac{1}{C_{\psi,\psi}}\int f(x)\overline{W_{\psi}\psi(y^{-1}
x)}\,dx.\nonumber
\end{align}
Hence $K(x,y)=W_{\psi}\psi(y^{-1}x)$ is the reproducing kernel for the 
image of $W_{\psi}$.

Moreover, \ref{eq1.2} shows $f(y)=C_{\psi,\psi}^{-1}\int f(x)(\psi,
\rho (x^{-1}y)\psi)\,dx=C_{\psi,\psi}^{-1}f\ast W_{\psi}(\psi)(y)$. 
Thus (g) and (h) follow.
\end{proof}

\section{Generalizations of the Wavelet Transform}\label{sec2}

Let $H$ be a locally compact Hausdorff topological group, and let 
$\pi:H\rightarrow\mathrm{GL}(n,\mathbb{R)}$ be a continuous 
homomorphism. Then we can define the semi-direct product $G=H\times_
{\pi}\mathbb{R}^{n}$ where the product is given by
\[
(a,x)(b,y)=(ab,x+\pi(a)y)\,.
\]

The group $G$ is locally compact; and if $dh$ is a left Haar measure on $H$
and $dx$ is the standard Lebesgue measure on $\mathbb{R}^{n}$, then a left
Haar measure on $G$ is given by
\[
d\mu_G(h,x)=\frac{\Delta_{\pi}(h)}{(2\pi)^{n}}\,dh\,dx
\]
where
\[
\Delta_{\pi}(h)=|\det(\pi(h))|^{-1}.
\]
Indeed, for continuous functions $f$ with compact support in $G$,
\begin{align*}
\int f((k,y)(h,x))\,d\mu_G (h,x)  &  =(2\pi)^{-n}\int\int f(kh,y+\pi
(k)x)\,|\det(\pi(k^{-1}kh))|^{-1}\,dh\,dx\\
&  =(2\pi)^{-n}\int\int f(h,y+\pi(k)x)\,|\det(\pi(k^{-1}h))|^{-1}\,dh\,dx\\
&  =(2\pi)^{-n}\iint f(h,y+x)\,|\det(\pi(h))^{-1}\det(\pi(k))|\,d(\pi
(k)^{-1}x)\,dh\\
&  =(2\pi)^{-n}\iint f(h,y+x)|\det(\pi(h))|^{-1}\,dx\,dh\\
&  =(2\pi)^{-n}\iint f(h,x)\,|\det(\pi(h))|^{-1}\,dx\,dh\\
&  =\int f(h,x)\,d\mu_G (h,x).
\end{align*}
As before, the action of $G$ on $\mathbb{R}^{n}$ given by
\[
T_{h,x}(y)=\pi(h)y+x
\]
induces a unitary representation $\rho$ of $G$ on $L^{2}(\mathbb{R}^{n})$.
Namely,
\[
\rho(h,x)f(y)=J(T_{h,x}^{-1})^{\frac{1}{2}}f(T_{h,x}^{-1}y)=\Delta_{\pi
}(h)^{1/2}f(\pi(h)^{-1}(y-x)).
\]

If $x\cdot y=\sum_{i=1}^{n}x_{i}y_{i}$ is the usual inner product on 
$\mathbb{R}^{n}$, then the Fourier transform defined formally by 

\[
\hat{f}(y)=\frac{1}{(2\pi)^{n/2}}\int f(x)e^{-ix\cdot y}\,dx
\]
is a unitary mapping from $L^{2}(\mathbb{R}^{n})$ onto itself.

To analyze the representation $\rho$, we again look at its Fourier 
transform $\hat{\rho}$. It is given by
\[
\hat{\rho}(h,x)\hat{f}=\widehat{\rho(h,x)f}.
\]
\begin{lemma}\label{rhohat2} Let $f\in L^2(\mathbb{R}^n)$. Denote by
$\pi(h)^t$ the transpose of the matrix $\pi(h)$
in $\mathrm{GL}(n,\mathbb{R})$. Then
$$\hat{\rho}(h,x)\hat{f}=\widehat{\rho(h,x)f}=
\Delta_{\pi}(h)^{-1/2}e^{-ix\cdot y}\,\hat{f}(\pi(h)^{t}y)
\, .$$
\end{lemma}
\begin{proof}
This is a simple calculation:
\begin{equation*}
\begin{split}
\hat{\rho}(h,x)\hat{f}(y)  &  =\frac{1}{(2\pi)^{n/2}}\int_{\mathbb{R}^{n}}
\rho(h,x)f(u)e^{-iu\cdot y}\,du\\
&  =\frac{1}{(2\pi)^{n/2}}\int_{\mathbb{R}^{n}}\Delta_{\pi}(h)^{1/2}
f(\pi(h^{-1})(u-x))e^{-iu\cdot y}\,du\\
&  =\frac{1}{(2\pi)^{n/2}}\int_{\mathbb{R}^{n}}\Delta_{\pi}(h)^{1/2}
f(\pi(h^{-1})u)e^{-i(u+x)\cdot y}\,du\\
&  =\frac{1}{(2\pi)^{n/2}}e^{-ix\cdot y}\int_{\mathbb{R}^{n}}\Delta_{\pi
}(h)^{1/2}f(\pi(h^{-1})u)e^{-iu\cdot y}\,du\\
&  =\frac{1}{(2\pi)^{n/2}}e^{-ix\cdot y}\int_{\mathbb{R}^{n}}\Delta_{\pi
}(h)^{1/2}f(u)e^{-i\pi(h)u\cdot y}\,d(\pi(h)u)\\
&  =\frac{1}{(2\pi)^{n/2}}e^{-ix\cdot y}\int_{\mathbb{R}^{n}}\Delta_{\pi
}(h)^{1/2}\,|\det(\pi(h))|\,f(u)e^{-iu\cdot\pi(h)^{t}y}\,du\\
&  =\Delta_{\pi}(h)^{-1/2}e^{-ix\cdot y}\,\hat{f}(\pi(h)^{t}y)\, .
\end{split}
\end{equation*}
\end{proof}

Let $A$ be a measurable subset of $\mathbb{R}^{n}$ with positive 
measure. Let $L_{A}^{2}(\mathbb{R}^{n})$ be the closed linear subspace 
of $L^{2} (\mathbb{R}^{n})$ consisting of those $L^{2}$ functions with 
$\hat{f}(y)=0$ for $y\notin A$. The subspaces $L_{A}^{2}(\mathbb{R}^ 
{n})$ are precisely the closed subspaces of $L^{2}(\mathbb{R}^{n})$ 
invariant under translation.

\begin{definition}
A $\pi(H)^{t}$ invariant measurable subset $A$ of $\mathbb{R}^{n}$ is 
ergodic (under Lebesgue measure) if any invariant measurable subset of 
$A$ has measure $0$ or has complement in $A$ with measure $0$.
\end{definition}

In particular, all homogeneous spaces obtained from the action of 
$\pi(H)^{t}$ are ergodic.

The action of $\pi(H)^{t}$ on $\mathbb{R}^{n}$ is said to have a {\bf 
measurable cross section} if there is a Borel set $E$ in $\mathbb{R}^ 
{n}$ meeting each set $\pi^{t}(H)y$ in precisely one point.

For completeness we present the following well known lemma.

\begin{lemma}\label{lem8} 
Suppose there is a Borel cross section $E$ for the action of $\pi(H)^
{t}$ on $\mathbb{R}^{n}$ and $A$ is an ergodic invariant Borel 
measurable set. Then $A$ contains a homogeneous space whose complement 
in $A$ has measure $0$.
\end{lemma}

\begin{proof}
Denote by $L$ the subgroup $\pi(H)^{t}$ of $\mathrm{GL}(n,\mathbb{R}) 
$. There is nothing to do if $A$ has measure $0$. Suppose $A$ has 
positive measure. Let $\mu$ be a probability measure equivalent to 
Lebesgue measure on $A$. Set $F=E\cap A$ where $E$ is a Borel set 
meeting each homogeneous space $Ly$ in exactly one point. Define $p: 
A\rightarrow F$ by $\{p(y)\}=Ly\cap F$. Then $p$ is a Borel function 
from $A$ onto $F$. Let $\mathcal{C}$ be a countable algebra of Borel 
subsets of $F$ which separate points in $F$. Set $\mathcal{C}_{0}=
\{W\in \mathcal{C}\mid\mu(p^{-1}(W))>0\}$. Since $p^{-1}(W)$ are 
invariant sets with positive measure, their complements in $A$ have 
measure $0$. Hence $\cap p^{-1}(W)$ has complement in $A$ with measure 
$0$. Thus $\cap p^{-1}(W)$ has positive measure. Suppose it contains 
two disjoint $L$-orbits $Ly$ and $Ly^{\prime}$. This implies $\cap W$ 
contains distinct points $p(y)$ and $p(y^{\prime})$. Since $\mathcal 
{C}$ is a separating family, there is a $U\in\mathcal{C}$ with $p(y)
\in U$ and $p(y^{\prime})\notin U$. Hence $U$ or $F-U$ belongs to 
$\mathcal{C}_{0}$. If $U\in\mathcal{C}_{0}$, then $p(y^{\prime})
\notin\cap W$ and if $F-U\in\mathcal{C}_{0}$, then $p(y)\notin\cap W$. 
Thus $\cap p^{-1}(W)$ is a homogeneous space whose complement in $A$ 
has measure $0$.
\end{proof}

\begin{theorem}\label{theorem9} Let $\pi : H\to \mathrm{GL}(n,\mathbb{R})$
be a continuous representation as before. Then
the following hold: 
\begin{itemize}
\item [(a)]\vskip -.3cm
A nonzero closed subspace $M$ of $L^{2}(\mathbb{R}^{n})$ is invariant 
under the representation $\rho$ if and only if $M=L_{A}^{2}(\mathbb 
{R}^{n})$ for some measurable $\pi(H)^{t}$ invariant subset $A$ of 
$\mathbb{R}^ {n}$ having positive measure. Moreover, this subspace is 
irreducible if and only if $A$ is ergodic.

\item[(b)] 
If the action of $\pi(H)^{t}$ on $\mathbb{R}^{n}$ has a Borel cross 
section, then every Borel measurable invariant ergodic subset $A$ of 
$\mathbb{R}^{n}$ contains a $\pi(H)^{t}$ homogeneous subset whose 
complement in $A$ has measure $0$. In particular, if $H$ is a Lie 
group, the irreducible subspaces of $\rho$ correspond to the open 
$\pi(H)^{t}$ homogeneous subspaces of $\mathbb{R}^{n}$.

\item[(c)] 
Let $A$ be ergodic under $\pi(H)^{t}$ with positive measure. Assume 
$\psi$ is a nonzero vector in $L_{A}^{2}(\mathbb{R}^{n})$. Then $\psi$ 
is a wavelet vector if and only if $h\mapsto\hat{\psi}(\pi(h)^{t}y)$ 
is in $L^{2}(H)$ for almost all $y\in A$. In particular, the 
stabilizer $H^{y}$ is compact for almost all $y\in A$.

\item[(d)] 
Assume $A$ is open and homogeneous under $\pi(H)^{t}$ and the 
stabilizer $H^{y}$ is compact for one (and hence all) $y\in A$. Then 
any nonzero $\psi$ with $\hat{\psi}\in C_{c}(A)$ is a wavelet vector
and the linear space of wavelet vectors in dense in $L^2_A(\mathbb{R}^n)$.
\end{itemize}
\end{theorem}

\begin{proof}
Let $M$ be a closed nonzero invariant subset of $L^2(\mathbb{R}^n)$. 
Set $S=\hat {M}$. By Lemma \ref{rhohat2}, $S$ is invariant under the unitary 
operators $\hat{\rho }(e,x)\hat {f}(y)=e^{-ix\cdot y}\hat {f}(y)$. Now 
every multiplication operator $M_hf=hf$ for $f\in L^2(\mathbb R^n)$ 
and $h\in L^{\infty}(\mathbb R^n)$ can be weakly approximated by a 
finite linear combination of the operators $\hat{\rho}(e,x)$. This 
implies $S$ is invariant under $M_h$ for all $h\in L^{\infty}(\mathbb 
R^n)$. In particular, if $P$ is the orthogonal projection onto $S$, 
then $PM_h=M_hP$ for all $h$. Hence $P=M_{\chi_A}$for some Borel 
measurable subset $A$ of $\mathbb{ R}^n$. Thus $M=L_A^2(\mathbb{R}^n) 
$.  

Since $S$ is invariant under $\rho$, we have $\hat{\rho}(h,0)(\chi_{A} 
f)=\chi_{A}\hat{\rho}(h,0)f$ for all $f\in L^{2}(\mathbb{R}^{n})$. 
Thus
\begin{align*}
\Delta_{\pi}(h)^{-\frac{1}{2}}\chi_{A}(\pi(h)^{t}y)f(\pi(h)^{t}(y))&  =
\hat{\rho}(h,0)(\chi_{A}f)(y)\\
&  =\chi_{A}(y)\hat{\rho}(h,0)f(y)\\
&  =\Delta_{\pi}(h)^{-\frac{1}{2}}\chi_{A}(y)f(\pi(h)^{t}y)
\end{align*}
a.e.\ $y$ for each $h\in H$. Thus $\chi_{A}(\pi(h)^{t}y)=\chi_{A}(y)$ 
a.e.\ $y$ for each $h$. Hence the sets $\pi(h^{-1})^{t}A$ and $A$ are 
essentially equal for all $h$.

Set $A^{\prime}=\{y\mid\pi(h)^{t}y\in A\mbox{\rm\ for a.e.\ }h\}$. By 
the left invariance of the measure $dh$, $A^{\prime}$ is invariant 
under $\pi^{t}(H)$. Moreover, by Fubini's Theorem, $A^{\prime}$ equals 
$A$ up to a set of measure $0$. Replacing $A$ by $A^{\prime}$, we see 
if $M$ is a nonzero closed $\rho$ invariant subspace of $L^{2} 
(\mathbb{R}^ {n})$, then $M=L_{A}^{2}(\mathbb{R}^{n})$ for some 
invariant Borel measurable $\pi(H)^{t}$ invariant subset of $\mathbb 
{R}^{n}$ with positive measure. The converse follows directly from 
Lemma \ref{rhohat2}. 

Next note the action of $\pi(H)^{t}$ on $A$ is not ergodic if and only 
if there exists an invariant subset $A_{0}$ of $A$ so that both $A_
{0}$ and $A-A_{0}$ have positive measure. This occurs if and only if 
$M$ has a proper invariant subspace of form $L_{A_{0}}^{2}(\mathbb{R}
^{n})$ which is equivalent to $M$ being reducible. We thus have (a).

The first statement in (b) follows from Lemma \ref{lem8}. To see the 
second, suppose $\pi(H)^{t}y$ is a homogeneous space with positive 
Lebesgue measure. Note $h\mapsto\pi(h)^{t}y$ will have range with 
positive measure in $\mathbb{\ R}^{n}$ if and only if the rank of this 
transformation is $n$. But this occurs only if $\pi(H)^{t}y$ is an 
open subset of $\mathbb{R}^{n}$.

To see (c) we set $W_{\psi}f=(f\mid\rho(h,x)\psi)$. By Lemma \ref{rhohat2}
\[
W_{\psi}f(h,x)=(\hat{f}\mid\hat{\rho}(h,x)\hat{\psi})=\int_{A}\hat{f}
(\omega)\Delta_{\pi}(h)^{-1/2}e^{ix\cdot\omega}\overline{\hat{\psi}(\pi
(h)^{t}\omega)}\,d\omega.
\]

Define $F(\omega)=\int_{H}|\hat{\psi}(\pi(a)^{t}\omega)|^{2}\,da$ for
$\omega\in A$. Note
\begin{align*}
F(\pi(h)^{t}\omega)  &  =\int_{H}|\hat{\psi}(\pi(a)^{t}\pi(h)^{t}\omega
)|^{2}\,da\\
&  =\int_{H}|\hat{\psi}(\pi(ha)^{t}\omega|^{2}\,da\\
&  =\int_{H}|\hat{\psi}(\pi(a)^{t}\omega|^{2}\,da\\
&  =F(\omega).
\end{align*}
The invariance of $F$ and the ergodicity of $A$ imply $F(\omega)$ is 
essentially constant on $A$; i.e., there is a constant $C_{\psi}^{2}
\leq\infty$ with
\[
\int_{H}|\hat{\psi}(\pi(a)^{t}\omega)|^{2}\,da=C_{\psi}^{2}\mbox{\rm\ for
a.e.\ }\omega\in A.
\]

Assume $f$ is in $L_{A}^{2}(\mathbb{R}^{n})$. Then 
\begin{align*}
\int|W_{\psi}f(h,x)|^{2}\,dx  &  =\int\int\int\hat{f}(\omega)\overline{\hat
{f}(\omega^{\prime})}\Delta_{\pi}(h)^{-1}\overline{\hat{\psi}(\pi(a)^{t}
\omega)}\hat{\psi}(\pi(h)^{t}\omega^{\prime})e^{ix\cdot\omega}e^{-ix\cdot
\omega^{\prime}}d\omega\,d\omega^{\prime}\,\,dx\\
&  =\Delta_{\pi}(h)^{-1}\int\int\int\hat{f}(\omega)\overline{\hat{\psi}
(\pi(a)^{t}\omega)}\,e^{ix\cdot\omega}\overline{\hat{f}(\omega^{\prime})}
\hat{\psi}(\pi(h)^{t}\omega^{\prime})e^{-ix\cdot\omega^{\prime}}
d\omega\,d\omega^{\prime}\,\,dx.
\end{align*}
Set $F(\omega)=\overline{\hat{f}(\omega)}\hat{\psi}(\pi(h)^{t}\omega)$. Then
\begin{align*}
\int|W_{\psi}f(h,x)|^{2}\,dx  &  =\Delta_{\pi}(h)^{-1}\int\int\int\hat
{f}(\omega)\overline{\hat{\psi}(\pi(a)^{t}\omega)}\,e^{ix\cdot\omega}
\overline{\hat{f}(\omega^{\prime})}\hat{\psi}(\pi(h)^{t}\omega^{\prime
})e^{-ix\cdot\omega^{\prime}}d\omega\,d\omega^{\prime}\,\,dx\\
&  =(2\pi)^{n}\Delta_{\pi}(h)^{-1}\int|\hat{F}(x)|^{2}\,dx\\
&  =(2\pi)^{n}\Delta_{\pi}(h)^{-1}\int_{A}|F(\omega)|^{2}\,d\omega\\
&  =(2\pi)^{n}\Delta_{\pi}(h)^{-1}\int_{A}|\hat{f}(\omega)|^{2}|\hat{\psi}
(\pi(h)^{t}\omega)|^{2}\,d\omega.
\end{align*}
Using Fubini's Theorem and
\[
\int_{H}|\hat{\psi}(\pi(h)^{t}\omega)|^{2}\,dh=C_{\psi}^{2}
\]
for a.e.\ $\omega$ in $A$, we see
\begin{align}
\int_{H\times_{\pi}\mathbb{R}^{n}}|W_{\psi}f(h,x)|^{2}\,d(h,x)  &  =\frac
{1}{(2\pi)^{n}}\int_{H}\Delta_{\pi}(h)\int|W_{\psi}f(h,x)|^{2}\,dx\,dh\nonumber\\
&  =\int_{H}\int_{A}|\hat{f}(\omega)|^{2}|\hat{\psi}(\pi(h)^{t}\omega
)|^{2}\,d\omega\,dh\nonumber\\
&  =\int_{A}|\hat{f}(\omega)|^{2}\int_{H}|\hat{\psi}(\pi(h)^{t}\omega
)|^{2}\,dh\,d\omega\label{fubini} \\
&  =C_{\psi}^{2}\int_{A}|\hat{f}(\omega)|^{2}\,d\omega\nonumber\\
&  =C_{\psi}^{2}(f,f).\nonumber
\end{align}
Hence $\psi$ is a wavelet vector if and only if $C_{\psi}^{2}<\infty$.

Now suppose $0<\int_H|\hat{\psi }(\pi (h)^ty)|^2\,dh<\infty$. Set $F(h)
=|\hat{\psi }(\pi (h)^ty)|^2$ and $K=H^y=\{a:\pi (a)^ty=y\}$. Note 
$F(ah)=F(a)$ for $a\in K$. Let $p:H\to K\backslash H$ be the natural 
projection and suppose $\mu$ is a regular right quasi-invariant 
measure on $K\backslash H$. Let $\int_{K\backslash H}m_{Kh}\,d\mu 
(Kh)$ be the disintegration of left Haar measure $m$ over the fibers of 
$p$.

Note  
\begin{align*}
m(E)&=m(aE)=\int_{K\backslash H}m_{Kh}(aE)\,d\mu (Kh)\\
&=\int_{K\backslash H}m_{Kh}(E)\,d\mu (Kh).
\end{align*}
Uniqueness of disintegrations yields $m_{Kh}(E)=m_{Kh}(aE)$ for all 
$E$ for each $a$ a.e.\ $Kh$. Hence if $\sigma :K\backslash H\to H$ is a 
Borel cross section, then since left Haar measures on $K$ are 
proportional, we have $$dm_{Kh}(a\sigma (Kh))=c(Kh)\,da$$ where $da$ 
is a left Haar measure on $K$ and $c(Kh)\ge 0$. Consequently, since 
$0<\int_HF(h)\,dh<\infty$ and $$\int F(h)\,dh=\int_{K\backslash H}c
(Kh)\int_KF(a\sigma (Kh))da\, d\mu (Kh)=\int_K\,da\,\int c(Kh)F(\sigma 
(Kh))\,d\mu (Kh),$$ we see Haar measure on $K$ is finite and thus $K$ 
is compact. We thus conclude $H^y$ is compact for almost all $y$.

Note (d) is consequence of $A$ is homeomorphic to $H/H^{y}$.
\end{proof}

\section{Wavelets and Symmetric Spaces}\label{sec3}
There are, as we will see in the next section, natural examples where 
$\mathbb{R}^n$ contains finitely many open orbits some of which do not 
have compact stabilizers. The obvious question is therefore, what do 
we do in these cases? The compactness of the stabilizer was used at one 
point in the proof of Theorem \ref{theorem9} to allow us to use 
Fubini's Theorem to interchange the integration over $A$ and $H$, see 
(\ref{fubini}): $$
\int_{H}\int_{A}|\hat{f}(\omega)|^{2}|\hat{\psi}(\pi(h)^{t}\omega
)|^{2}\,d\omega\,dh=
\int_{A}|\hat{f}(\omega)|^{2}\int_{H}|\hat{\psi}(\pi(h)^{t}\omega
)|^{2}\,dh\,d\omega\, .$$ We notice that the right hand side makes 
perfectly sense if $h\mapsto |\hat{\psi}(\pi(h)^{t}\omega )|^{2}$ is 
integrable over $H/H^{\omega}$, where $H^{\omega}=\{h\in H\mid \pi(h) 
^t\omega =\omega \}$. The problem is, that the stabilizer group 
$H^\omega$ depends on the point $\omega$. We will show how to overcome 
this difficulty in case where $H^\omega$ is a \textit{symmetric} 
subgroup of $H$. We will show that in this case we can replace $ H$ by 
a closed subgroup $R$ such that the open orbit essentially decomposes 
into finitely many smaller open orbits which are homogeneous under $R$ 
and have compact stabilizers. The drawback is that in most cases the 
group $R$ depends on the orbit, but at least this allows us to 
reconstruct the function from its wavelet transform. The tools that we 
use are the structure theory of reductive symmetric spaces 
corresponding to an involution $\tau : H\to H$. In particular we will 
need a \textit{Cartan involution} $\theta : H\to H$ commuting with 
$\tau$ and the structure of minimal $\theta\tau$-stable parabolic 
subgroups. To avoid technical details we will only discuss this for a 
linear Lie group invariant under transposition. Then $h\mapsto (h^ 
{-1})^t$ is a Cartan involution on $H$. We refer to \cite{GO87} for 
the discussion of the general case. We thus will be dealing with 
closed subgroups $H$ of $\mathrm{GL}(n,\mathbb{R})$ invariant under 
transposition. A conjugate being invariant under transposition is 
one of several alternative definitions for reductive linear groups.

\begin{definition}
A closed subgroup $H\subset\mathrm{GL}(n,\mathbb{R})$ is called
\textbf{reductive} if there exists a $x\in\mathrm{GL}(n,\mathbb{R})$ 
such that $xHx^{-1}$ is invariant under transposition, $a\mapsto a^
{t}$.
\end{definition}

We will now assume that $H$ is reductive. For simplicity we can then 
assume that $H^t=H$. Let $\pi$ be the natural action of $H$ on 
$\mathbb{R}^n$ given by $\pi (h)v=hv$. If $\rho$ is the natural 
representation of $H\times_{\pi}\mathbb R^n$ on $L^2(\mathbb R^n)$, we 
will be interested in the irreducible subrepresentations of $\rho$ 
corresponding to open orbits in $\mathbb R^ n$ under the 
contragredient action $h\cdot v=(h^{-1})^tv$. We recall if $\mathcal 
O$ is such an orbit, then the subspace $L^2_{\mathcal O}(\mathbb R^n)
=\{f\in L^2(\mathbb R^n)|\,\hat {f}\text{\rm \ vanishes off }
\mathcal O\}$ is an irreducible subspace.  In order to handle cases 
where $\mathcal O=H\cdot u$ do not have compact stabilizers, we will 
assume that there exists an involution $\tau :H\rightarrow H$ such 
that 
\[
H_{o}^{\tau}\subset H^{u}=\left\{  h\in H\mid h\cdot u=u\right\}  \subset
H^{\tau}
\]
where
\[
H^{\tau}=\left\{  h\in H\mid\tau(h)=h\right\}
\]
and the subscript $_{o}$ indicates the connected component containing 
the identity element. Thus $\mathcal{O}=H/H^{u}$ is \textit{a 
semisimple symmetric space}. We can assume that $L=H^{u}$ is also 
invariant under transposition. Then $\theta:h\mapsto(h^{t})^{-1}$ and 
$\tau$ commute. Let $K=\mathrm{O} (n)\cap H=\left\{ k\in H\mid k^{t} 
=k^{-1}\right\} $. Then $K$ is a maximal compact subgroup of $H$ and 
$L\cap K$ is a maximal compact subgroup of $L$. Denote by $\frak{h}$ 
the Lie algebra of $H$. One has
\[
\frak{h}=\left\{  X\in M_{n}(\mathbb{R})\mid 
\forall t\in\mathbb{R}:e^{tX}\in H \right\}.
\]
Then $\frak{h}$ decomposes as
\begin{align*}
\frak{h}  &  =\frak{k}\oplus\frak{s}\\
&  =\frak{l}\oplus\frak{q}\\
&  =\frak{k}\cap\frak{l}\oplus\frak{k}\cap\frak{q}\oplus\frak{s}\cap
\frak{l}\oplus\frak{s}\cap\frak{q}\,
\end{align*}
where
$$\frak{s}=\left\{ X\in\frak{h}\mid X^{t}=X\right\}$$
is the 
subspace of symmetric matrices, and
$$\frak{q}=\left\{ X\in\frak{h}
\mid\tau(X)=-X\right\}\, . $$
Notice that
\[
\lbrack\frak{l},\frak{q}]\subset\frak{q\quad}\text{and}\quad\lbrack
\frak{k},\frak{s}]\subset\frak{s\,.}
\]
Recall the linear maps $\mathrm{Ad}(a),\mathrm{ad}(X):\frak{h}
\rightarrow\frak{h}$, $a\in H$, $X\in\frak{h}$, are given by 
$\mathrm{ad}(X)Y=XY-YX$ and $\mathrm{Ad}(a)Y=aYa^{-1}$. Let $\frak{a}$ 
be a maximal commutative subspace of $\frak{s}\cap\frak{q}$; thus 
$XY-YX=0$ for all $X,Y\in\frak{a}$. Then the algebra $\mathrm{ad} 
(\frak{a})$ is also commutative. Define an inner product 
$(\cdot\mid\cdot)$ on $\frak{h}$ by $(X\mid Y): =\mathrm{Tr}(XY^{t})$. 
Then
\begin{align*}
(\mathrm{ad}(X)Y\mid Z)  &  =\mathrm{Tr}(XY-YX)Z^{t}\\
&  =-\mathrm{Tr}(Y(ZX^{t}-X^{t}Z)^{t})\\
&  =(Y\mid\mathrm{ad}(X^{t})Z)\,.
\end{align*}
Thus $\mathrm{ad}(X)^{t}=\mathrm{ad}(X^{t})$. In particular, if 
$X\in\frak{s}$ then $\mathrm{ad}(X)$ is symmetric. It follows that we 
can diagonalize the action of $\frak{a}$ on $\frak{h}$. Specifically, 
for $\alpha\in\frak{h}^{\ast}$ set
\[
\frak{h}_{\alpha}=\left\{  Y\in\frak{h}\mid\forall X\in\frak{a}\,:\mathrm{ad}
(X)Y=\alpha(X)Y\right\}  \,.
\]
Let $\Delta=\left\{\alpha\in\frak{h}^{\ast}\mid\frak{h}_{\alpha}
\not =\left\{0\right\}\right\}\setminus\left\{0\right\}$. Notice
that the set $\Delta$ is finite. Hence there is a $X_r\in\frak{a}$ 
such that $\alpha(X_r)\not =0$ for all $\alpha\in\Delta$. Let 
$\Delta^{+}=\left\{\alpha\mid\alpha(X_r)>0\right\}$ and 
$\Delta^{-}=\left\{\alpha\mid\alpha(X_r)<0\right\}=-\Delta^{+}$. Let
\[
\frak{n}=\bigoplus_{\alpha\in\Delta^{+}}\frak{h}_{\alpha}
\quad\mathrm{and\quad}\frak{\bar{n}}=
\bigoplus_{\alpha\in\Delta^{-}}\frak{h}_{\alpha}\,.\quad
\]
Let $\frak{m}_{1}=\left\{  X\in\frak{h}\mid\lbrack\frak{a},X]=\left\{
0\right\}  \right\}  $, and $\frak{m}=\left\{  X\in\frak{m}_{1}\mid\forall
Y\in\frak{a}:(X\mid Y)=0\right\}  $. Then $\frak{m}_{1}=\frak{m}\oplus
\frak{a}$. Furthermore $\frak{\bar{n}}$, $\frak{m}$, $\frak{n}$, and
$\frak{p}=\frak{m}\oplus\frak{a}\oplus\frak{n}$ are subalgebras of $\frak{h}$
and
\begin{align*}
\frak{h}  &  =\frak{\bar{n}}\oplus\frak{m}\oplus\frak{a}\oplus\frak{n}\\
&  =\frak{k}+\frak{p}.
\end{align*}
Notice that the last sums are not direct because $\frak{k}\cap\frak{p} 
=\frak{k}\cap\frak{m}$ and $\frak{l}\cap\frak{p} =\frak{l}\cap\frak
{m}$. Furthermore $P\cap L$ is not necessarily compact. To deal with 
this, we let $\frak{m}_{2}$ be the smallest subalgebra containing 
$\frak{m}\cap\frak{s}$. Thus $\frak{m}_{2}=[\frak{m}\cap\frak{s},
\frak{m}\cap\frak{s}]\oplus\frak{m}\cap\frak{s}$.

\begin{lemma}
$\frak{m}_{2}$ is an ideal in $\frak{m}$ and contained in $\frak{m}
\cap\frak{l}.$
\end{lemma}

\begin{proof}
That $\frak{m}\cap\frak{s}\subset\frak{l}$ $\ $follows from the fact 
that $\frak{a}$ is maximal abelian in $\frak{q}\cap\frak{s}$. Hence 
$\frak{m} _{2}\subset\frak{l}$ follows from the fact that $\frak{l}$ 
is an algebra. By definition we have $\frak{m}=\frak{m}\cap\frak{k}
+\frak{m}_{2}$. We have by definition $[\frak{m}_{2},\frak{m}_{2}]
\subset\frak{m}_{2}$. Therefore we only have to show that $[\frak{m}
\cap\frak{k},\frak{m}_{2}]\subset\frak{m}_{2}$. But $[\frak{k},\frak
{s}]\subset\frak{s}$ and hence $[\frak{m\cap k} ,\frak{m}\cap\frak{s}
]\subset\frak{m}\cap\frak{s}$. The claim follows now by the Jacobi 
identity and the fact that $\frak{m}_{2}$ is generated by $\frak{m\cap 
s}$.
\end{proof}

Let
\[
\frak{b}:=\left\{  X\in\frak{m}\mid\forall Y\in\frak{m}_{2}:(X,Y)=0\right\}
=\frak{m}_{2}^{\perp}.
\]
Then $\frak{b}$ is an ideal in $\frak{m}$ and $\frak{m}=\frak{b}\oplus
\frak{m}_{2}$. Let
$$N_{K}(\frak{a})=\left\{  k\in K\mid\forall X\in
\frak{a}:\mathrm{Ad}(a)X\in\frak{a}\right\}  $$
and
$$H_{L\cap K}(\frak{a}
)=N_{K}(\frak{a})\cap L\, .$$
Finally let
$$M_{K}=Z_{K}(\frak{a})=\left\{ 
k\in K\mid\forall X\in\frak{a}:aXa^{-1}=X\right\}$$
and
$$M_{H}=Z_{H}
(\frak{a})\, .$$ Then
\begin{equation}\label{Weyl}
W=N_{K}(\frak{a})/M_{K}
\end{equation}
is a finite group. Let 
\begin{equation}\label{smallWeyl}
W_{0}=N_{K\cap L}(\frak{a})/M_{K}\cap L\, .
\end{equation}
Then $W_0$ is a subgroup of $W$. Choose $w_{0}=e$, $w_
{1},\ldots,w_{r}\in W$ such that
\begin{equation}\label{represent}
W=\bigcup w_{j}W_{0}\quad\text{(disjoint union)\thinspace.}
\end{equation}
Let $s_{j}\in N_{K}(\frak{a})$ be such that $s_{j}M_{K}=w_{j}$, and 
$s_{0}=e$. Let $P=\left\{ a\in H\mid\mathrm{Ad}(a)\frak{p}=\frak{p}
\right\} $. Then $P$ is a closed subgroup of $H$. Let $A=\left\{ e^
{X}\mid X\in\frak{a}\text{ }\right\} $ and $N=\left\{ e^{X}\mid 
X\in\frak{n}\right\} $. Then $A$, and $N$ are closed subgroups of $P$. 
Let $M_{2}$ be the group generated by $\exp(\frak{m}_{2})$, and $B_
{o}=\exp(\frak{b})$. Then $F=\exp(i\frak{a})\cap K\subset M_{K}$ is 
finite and such that $B=FB_{o}$ is a group. Furthermore
\[
B\times M_{2}\times A\ni(b,m,a)\mapsto bma\in Z_{H}(A)
\]
is a diffeomorphism. Notice that by definition $F$ is central in 
$Z_{H}(A)$ and $mFm^{-1}=F$ for all $m\in N_{K}(\frak{a})$. Let 
$M=BM_{2}$. Then $P=MAN$. Furthermore each element $p\in P$ has a unique 
expression $p=man$ with $m\in M$, $a\in A$, $n\in N$. The group $P$ is 
called a \textit{minimal }$\theta\tau$ \textit{stable parabolic 
subgroup of }$H$. This is still not the correct group for us to work 
with because $P\cap s_{j}Ls_{j}^{-1}$ may not be compact. We 
therefore set
\[
R=BAN\,.
\]
Then $R$ is a closed subgroup of $H$ with Lie algebra $\frak{r}=\frak{b}
\oplus\frak{a}\oplus\frak{n}$. Notice that $\frak{h}=\frak{l}+
(\frak{b}\oplus\frak{a}\oplus\frak{n})$ and 
$\frak{l}\cap(\frak{b}\oplus\frak{a}\oplus\frak{n})
=\frak{l}\cap\frak{b}$.

\begin{lemma}
Let the notation be as above. Then the following hold:
\begin{enumerate}
\item\vskip -.2cm  
If $s\in N_{K}(\frak{a})$, then $sBs^{-1}=B$ and $sM_{2}s^{-1}=M_{2}$.
\item 
$P\cap s_{j}Ls_{j}^{-1}=M\cap s_{j}Ls_{j}^{-1}$.
\item  
We have $R\cap s_{j}Ls_{j}^{-1}\subset K$ is compact.
\end{enumerate}
\end{lemma}
\begin{proof}
\ \linebreak
(a) First we notice that if $X\in\frak{m}$, $Y\in\frak{a}$ and $s\in
N_{K}(\frak{a})$ then
\[
\lbrack\mathrm{Ad}(s)X,Y]=\mathrm{Ad}(s)[X,\mathrm{Ad}(s)^{-1}Y]=0\,.
\]
Hence $\mathrm{Ad}(s)\frak{m}=\frak{m}$. Furthermore $\mathrm{Ad}
(s)\frak{s}=\frak{s}$ as $s\in K$. It follows that $\mathrm{Ad}(s)\frak{m}
_{2}=\frak{m}_{2}$. As $\mathrm{Ad}(s)$ is an orthogonal transformation it
follows that $\mathrm{Ad}(s)\frak{b}=\frak{b}$. Thus $sB_{o}s^{-1}=B$ and
$sM_{2}s^{-1}=M_{2}$. As $sFs^{-1}=F$ it follows that $sBs^{-1}=B$.

(b) See \cite{GO87}.

(c) Let $p\in R\cap s_{j}Ls_{j}^{-1}$. Then again by using
\cite{GO87} one has $p\in M\cap L\cap B\subset M\cap K$. 
\end{proof}

We are now able to state the following classical result of Matsuki
\cite{Ma79}.

\begin{theorem}
[Matsuki]The sets $Ps_{0}L,\ldots,Ps_{r}L$ are disjoint and open. 
Furthermore $\bigcup_{j=0}^{r}Ps_{j}L$ is open and dense in $H$.
\end{theorem}

As a consequence of this we get (cf. \cite{GO87}, p. 611):
\begin{corollary}
The sets $Rs_{0}L,\ldots,Rs_{r}L$ are disjoint and open. Furthermore 
$\bigcup_{j=0}^{r}Rs_{j}L$ is open and dense in $H$.
\end{corollary}

\begin{proof}
This follows from the fact that $Ms_{j}L=Bs_{j}L$.
\end{proof}
One can actually state more than this, see \cite{GO87}.

The groups $H$ and $L$ being reductive are unimodular. This fact
is important for the integral formulas in the next theorem.
\begin{theorem} Let the notation be as above. Then the following
hold:
\begin{enumerate}
\item \vskip -.2cm 
There exists a real analytic function $p:H\rightarrow\mathbb{R}_{0} ^
{+}=\left\{ t\in\mathbb{R}\mid t\geq0\right\} $ such that
\[
\bigcup_{j=0}^{r}Rs_{j}L=\left\{  x\in H\mid p(x)>0\right\}  \,.
\]
In particular it follows that $H\setminus\bigcup_{j=0}^{r}Rs_{j}L$ has 
measure zero in $H$.

\item  
We can normalize the Haar measures on the groups $R$ and $L$ so that
\[
\int_{H}f(x)\,dh=\int_{R}\int_{L}f(rs_{j}l)\,dldr
\]
for all $f\in L^{1}(H)$ with $\mathrm{Supp}(f)\subset Rs_{j}L$.
\end{enumerate}
\end{theorem}

{}From this we now get the following result:

\begin{theorem}
Let $H$ be a reductive group acting on $\mathbb{R}^{n}$. Assume that 
$\mathcal{O}=H\cdot u$ is an open orbit and $H^{u}$ is noncompact and 
symmetric. Then there exists a closed subgroup $R$ of $H$ and elements 
$y_{0}=y,\ldots,y_{r}\in\mathcal{O}$ such that the following hold:
\begin{enumerate}
\item  \vskip -.2cm
The orbits $\mathcal{U}_{j}=R\cdot y_{j}\subset\mathcal{O}$ are open and 
disjoint;
\item  
If $x\in\mathcal{U}_{j}$, then $R^{x}=\left\{ r\in R\mid r\cdot x=x\right\}$ 
is compact;

\item 
$\bigcup\mathcal{U}_{j}\subset\mathcal{O}$ is open and dense, and the 
complement has measure zero;

\item  
The space $L_{\mathcal{O}}^{2}(\mathbb{R}^{n})$ decomposes 
orthogonally as a representation of $R\times_{\pi}\mathbb{R}^n$ 
into irreducible parts
\[
L_{\mathcal{O}}^{2}(\mathbb{R}^{n})=\bigoplus_{j=0}^{r}L_{\mathcal{U}_{j}}
^{2}(\mathbb{R}^{n})\,.
\]

\item  
The space of $R\times_{\pi}\mathbb{R}^n$-wavelets is dense in each of 
the spaces $L_{\mathcal {U}_{j}}^{2}(\mathbb{R}^{n})$.
\end{enumerate}
\end{theorem}

Thus $\rho$ on $L^2_{\mathcal O}(\mathbb R^n)$ restricted to 
$G_0=R\times_ {\pi} \mathbb R^ n$ decomposes into a finite sum 
irreducible representations, one for each $\mathcal{U}_j$.  

Let $P_j$ be the orthogonal projection of $L^2(\mathbb R^n)$ onto 
$L^2_{\mathcal U_j}(\mathbb R^n)$. Then $\rho_j=\rho |_{G_0}$ on 
$L^2_{\mathcal U_j}(\mathbb R^n)$ is irreducible and any nonzero 
$\psi_j$ with $\hat{\psi}_j\in C_c(\mathcal U_j)$ is a $G_0$ 
wavelet for $\rho_ j$.  

Thus the mapping
$$W_j:L^2_{\mathcal U_j}(\mathbb R^n)\to L^2(G_0)$$
defined by
$$W_jf(g)=(f\mid \rho (g)\psi_j)$$
is a intertwining operator of $\rho_j$ into $L^2(G_0)$. Moreover, there
is a scalar $C_j$ such that
\[
(W_{j}f\mid W_{j}g)=C_{j}^{2}(f\mid g)
\]
for $f,g\in L^2_{\mathcal U_j}(\mathbb R^n)$.

Also for $f\in L^{2}(\mathbb{R}^{n})$,
\[
P_{j}f=\frac1{C_{j}^{2}}\rho_{0}(W_{j}f) =\frac1{C_{j}^{2}}\int_{G_{0}}
(f\mid \rho_{0}(g)\psi_{j})\,\rho_{0}(g)\psi_{j}\,dg.
\]
Hence if $f\in L^2_{\mathcal O}(\mathbb R^n)$, we have
\[
f=\sum P_{j}f=\sum_{j=0}^{r}\frac1{C_{j}^{2}} \int_{G_{0}}(f\mid 
\rho(g)\psi_{j})\rho(g)\psi_{j}\,dg.
\]

By normalizing the $\psi_j$, we may assume $C_j=1$ for all $j$. Then 
if $f\in L^2(\mathbb{R}^n)$ and $P=\sum P_j$, then 
\[Pf=\int_{G_{0}}\sum_{j=0}^{r}(f\mid \rho(g)\psi_{i})\rho(g)\psi_{j}\,dg.
\]We thus can reconstruct $Pf$ weakly from the finite collection of 
wavelet transforms $W_j(f)$ which we call a wavelet package. This 
remains possible even if there are no genuine wavelets for $G$ on the 
orbit $\mathcal{O}$. 

\begin{remark}
The authors thought it feasible to describe which functions are 
wavelets in the symmetric case by determining if their Fourier 
transforms were square integrable over the homogeneous orbits 
$H^u\backslash H$. However, the integrals formulas $\int_ R|\hat{\psi 
}(\pi (r)^ty_j)|^2\,dr$ involve left Haar measure and the measure 
needed on $H^u\backslash H$ involve the right invariant measure on 
$R$. The nonunimodularity of $R$ prevents these from being matched up.
\end{remark}

\section{Examples}\label{sec4}

In this section we present two often cited examples.
We discuss the first example of the group
$H=\mathbb{R}^+\mathrm{SO}_o(1,n)$ acting
on $\mathbb{R}^{n+1}$ in details.
For a very clear and detailed discussion of
cases where our assumptions are satisfied, i.e., the
pre-homogeneous vector spaces,  parabolic case, and
the explicite construction of the involutions corresponding
to the stabilizer of  points in an open orbit, we refer
to the preprint by N. Bopp and H. Rubenthaler \cite{BR}.
The varyity of those examles shows,
that the symmetric case gives a class of wavelets,
and wavelet transforms, worth studying in more details.

We start by  taking the group $H$ to be $\mathbb{R}^{+}\mathrm{SO}_o(1,n)$ and let 
it act on $\mathbb{R}^{1+n}$ in its usual way. For the second we 
let $H$ be the group $\mathrm{GL}(n,\mathbb{R})$ and let it act on the 
Euclidean space $\mathrm{Symm}(n,\mathbb{R})$ of symmetric $n\times 
n$-matrices by $a\cdot X:=aXa^t$.  

We mention there are other further examples where $H$ is an 
automorphism group of a symmetric open convex cone; e.g., take $H$ to 
be $\mathbb{R}^{+}$\textrm{SL}$(n,\mathbb{C})$ acting on the space 
$\mathrm{Symm}(n,\mathbb{C})$ of symmetric $n\times n$-matrices by 
$a\cdot X=aXa^{\ast}$, or by taking $H$ to be $\mathbb{R}^{+}\mathrm{ 
SU}^{\ast}(2n)$ or $\mathbb{R}^{+}E_{6(-26)}$ with similarly defined 
actions.

\vskip 5pt
{\sc\bf\large Example} ${\mathbf (\mathbb{R}^{+}\mathrm{SO}_{o}(1,n),
\mathbb{R}^{1+n})}$:
\vskip 3pt

Let $H=\mathbb{R}^{+}\mathrm{SO}(1,n)$ where $\mathbb{R}^{+}$ stands 
for the group of matrices $\{\lambda\mathrm{I}_{1+n}:\lambda>0\}$ and 
$\mathrm{I}_ {k}$ is the $k\times k$ identity matrix. Notice that 
$\mathbb{R}^{+}$ is central in $H$. The group $H$ acts on $\mathbb 
R^{1+n}$ by matrix multiplication on $\mathbb{R}^{1+n}$. If 
$\lambda\in\mathbb{R}^{+}$ and $h\in\mathrm {SO}(1,n)$, then  
\[
\Delta(\lambda h)=\det(\lambda h)^{-1}=\lambda^{-(1+n)}.
\]

We write elements in $\mathbb{R}^{1+n}$ as $(t,x)$ with $t\in\mathbb{R}$
and $x\in\mathbb{R}^{n}$. If one defines a bilinear form $\beta$ on
$\mathbb {R}^{1+n}$ by
\[
\beta((t,x),(s,y))=ts-(x\mid y)
\]
where $(\cdot\mid\cdot)$ is the usual inner product on $\mathbb{R}^
{n}$, then $\mathrm{SO}_{o}(1,n)$ is the connected component of the
identity of the group
\[
\mathrm{O}(\beta):=\left\{  g\in\mathrm{GL}(1+n,\mathbb{R})\mid\forall
u,v\in\mathbb{R}^{1+n}:\beta(gu,gv)=\beta(u,v)\right\}  \,.
\]

The Lie algebra of $\mathrm{SO}_o(1,n)$ consists of all matrices $n+1$
by $n+1$ matrices $X$ satisfying $X^tB+BX=0$ where $B=\left(\begin{matrix}
1&0\\ 0&-I_n\end{matrix} \right)$.

We will when convenient write elements in $\mathrm{SO}_{o}(1,n)$ in
block form
\[
g=\left(
\begin{array}
[c]{cc}
a & x\\
y^{t} & A
\end{array}
\right)  \,,\quad a\in\mathbb{R},x,y\in\mathbb{R}^{n},A\in M_{n}(\mathbb{R})
\]
Notice that actually $a>0$ and that $H^{t}=H$.

Denote the matrix $\left(\delta_{\nu,i}\delta_{\mu,j}\right)_{i,j=1}^
{n+1}$ by $E_{\nu,\mu}$ and let $X=E_{1,n+1}+E_{n+1,1}\in\frak{s}$.
Let $\frak{a}=\mathbb{R}X+\mathbb{R}\mathrm{I}_{n+1}$, where the
second factor corresponds to the dilations. Then $\frak{a}$ is a
maximal abelian subalgebra of $\frak{s}$. We will write
\[
a_{t}=e^{tX}=\left(
\begin{array}
[c]{ccc}
\cosh(t) & 0 & \sinh(t)\\
0 & \mathrm{I}_{n-1} & 0\\
\sinh(t) & 0 & \cosh(t)
\end{array}
\right)  \in A\,.
\]
 Imbed $\mathrm{SO}_{o}(n)$ into $\mathrm{GL}(1+n,\mathbb{R})$ by
\[
k\mapsto\left(
\begin{array}
[c]{cc}
1 & 0\\
0 & k
\end{array}
\right)  \,.
\]
We will also imbed $\mathrm{SO}_{o}(1,n-1)$ into $\mathrm{SO}_{o}(1,n)$ by
\[
g^{\prime}\mapsto g=\left(
\begin{array}
[c]{cc}
g^{\prime} & 0\\
0 & 1
\end{array}
\right)  \,.
\]
Then $\mathrm{SO}_{o}(n)\subset\mathrm{SO}_{o}(1,n)$ and $K=\mathrm{SO}
_{o}(n)$ is a maximal compact subgroup of $H$. Let
\begin{align*}
\mathcal{O}_{1}  &  =\left\{  u=(t,x)\in\mathbb{R}^{1+n}\mid\beta
(u,u)>0,t>0\right\} \\
\mathcal{O}_{2}  &  =\left\{  u=(t,x)\in\mathbb{R}^{1+n}\mid\beta
(u,u)>0,t<0\right\} \\
\mathcal{O}_{3}  &  =\left\{  u=(t,x)\in\mathbb{R}^{1+n}\mid\beta
(u,u)<0\right\} \\
\mathcal{C}  &  =\left\{  u=(t,x)\in\mathbb{R}^{1+n}\mid\beta(u,u)=0\right\}
\end{align*}
Then $\mathcal{O}_{1}$, $\mathcal{O}_{2}$, and $\mathcal{O}_{3}$ are open,
and the complementary set $\mathcal{C}$ is the zero set of $\beta$ and
is a union of orbits of smaller dimension. We have the
decomposition:
\[
\mathbb{R}^{1+n}=\mathcal{O}_{1}\cup\mathcal{O}_{2}\cup\mathcal{O}_{3}
\cup\mathcal{C\quad}\text{(disjoint union)\thinspace.}
\]

\begin{lemma}
The sets $\mathcal{O}_{1}$, $\mathcal{O}_{2}$, and $\mathcal{O}_{3}$ are
homogeneous under $H$. Furthermore the following hold:
\begin{enumerate}
\item \vskip -.2cm 
The sets $\mathcal{O}_{1}$ and $\mathcal{O}_{2}$ are homogeneous self
dual convex cones.

\item  
If $u\in\mathcal{O}_{1}\cup\mathcal{O}_{2}$, then $H^{u}\simeq
\mathrm{SO}_{o}(n)$ is compact.

\item  
If $u\in\mathcal{O}_{3}$, then $H^{u}\simeq\mathrm{SO}_{o}(1,n-1)$ is
noncompact.

\item  
If we replace $H$ by the non-connected group $H^{\ast}=\mathbb{R} ^
{\ast}\mathrm{SO}_{o}(1,n-1)$, then $\mathcal{U}_{2}=\mathcal{O}_{1}
\cup\mathcal{O}_{2}$ is homogeneous.
\end{enumerate}
\end{lemma}

\begin{proof}
All of this is well known so that we only prove that $\mathcal{O}_{3}$
is homogeneous and that $H^{u}\simeq\mathrm{SO}_{o}(1,n-1)$ for $u\in
\mathcal{O}_{3}$. Considering vectors in $\mathbb{R}^k$ as column vectors,
we let $u=e_{n+1}=(0,\ldots,0,1)^{t}$. In the following we will also
view $u$ as a vector in $\mathbb{R}^{n}\subset\mathbb{R} ^{1+n}$. Then
$u\in\mathcal{O}_{3}$. Notice that
\[
a_{t}u=(\sinh(t),0,\ldots,0,\cosh(t))^{t}\,.
\]
Let $v=(s,y)\in\mathcal{O}_{3}$. By multiplying by $\lambda=\left(
-\beta(v,v)\right) ^{-1/2}$, we may assume that $\beta(v,v)=-1$. Then
$||y||^{2}-s^{2}=1$. Hence we can find $t$ such that
\[
||y||=\cosh(t)\qquad\text{and}\qquad s=\sinh(t).
\]
In particular, both $\cosh(t)u$ and $y$ are in $S_{\cosh(t)}(0)=
\left\{x\in\mathbb{R}^{n}\mid||x||=\cosh(t)\right\}
\subset\mathbb{R}^{1+n}$. 
Recall that $\mathrm{SO}_{o}(n)$ acts transitively on $S_{\cosh(t)}
(0)$. Choose $k\in\mathrm{SO}_{o}(n)\subset\mathrm{SO}_{o}(1,n)$ with
$k(\cosh(t)u)=y$. Then $ka_{t}u=v$. Thus $\mathrm{SO}_{o}(1,n)$ acts
transitively. 

Let $g=\lambda\left(
\begin{array}
[c]{cc}
a & x^{t}\\
y & A
\end{array}
\right)\in H$ satisfy $gu=u$. Note first that $\lambda=1$ and thus
\[
gu=\left(
\begin{array}
[c]{c}
x_{n}\\
a_{1n}\\
\vdots\\
a_{nn}
\end{array}
\right)  =\left(
\begin{array}
[c]{c}
0\\
0\\
\vdots\\
1
\end{array}
\right)
\]
Hence $x\in\mathbb{R}^{n-1}$ and $A$ has the form
\[
A=\left(
\begin{array}
[c]{cc}
A^{\prime} & 0\\
0 & 1
\end{array}
\right)  \,.
\]
That the last row is $(0,\ldots,0,1)$ follows from the fact that
$g\in\mathrm{SO}(1,n)$. Since $g\in\mathrm{SO}_{o}(1,n)$, one has
\[
\left(
\begin{array}
[c]{cc}
a & 0\\
0 & A^{\prime}
\end{array}
\right)  =g^{\prime}\in\mathrm{SO}_{o}(1,n-1)\subset\mathrm{SO}_{o}(1,n)\,.
\]
Hence the claim.
\end{proof}

>From this it follows, that we can decompose $L^{2}(\mathbb{R}^{1+n})$
as a representation of $G=H\times_{\pi}\mathbb{R}^{n+1}$ into irreducible
parts as
\[
L^{2}(\mathbb{R}^{n})\simeq_{G}L_{\mathcal{O}_{1}}^{1}(\mathbb{R}^{n})\oplus
L_{\mathcal{O}_{2}}^{1}(\mathbb{R}^{n})\oplus L_{\mathcal{O}_{3}}
^{1}(\mathbb{R}^{n}).
\]
The first two spaces have $G$-wavelets. If we replace $H$ by 
$\mathbb{R}^*\mathrm{SO}_o(n+1)$, then this decomposition simplifies to
\[
L^{2}(\mathbb{R}^{n})\simeq_{G}L_{\mathcal{O}_{1}\cup\mathcal{O}_{2}
}^{1}(\mathbb{R}^{n})\oplus L_{\mathcal{O}_{3}}^{1}(\mathbb{R}^{n})
\]
where the first space has a dense subspace of wavelets. To deal with
the orbit $\mathcal{O}_{3}$ we show that $H^{e_{n+1}}$ is in fact a
symmetric subgroup of $H$. To that end define an involution $\tau$ on $H$
by
\begin{align}
\tau\left(  \lambda\left(
\begin{array}
[c]{cc}
A & x\\
v^{t} & a
\end{array}
\right)  \right)   &  =\frac{1}{\lambda}\left(
\begin{array}
[c]{cc}
\mathrm{I}_{n} & 0\\
0 & -1
\end{array}
\right)  \left(
\begin{array}
[c]{cc}
A & x\\
v^{t} & a
\end{array}
\right)  \left(
\begin{array}
[c]{cc}
\mathrm{I}_{n} & 0\\
0 & -1
\end{array}
\right) \nonumber\\
&  =\frac{1}{\lambda}\left(
\begin{array}
[c]{cc}
A & -x\\
-v^{t} & a
\end{array}
\right)  \,. \label{eqtau}
\end{align}
Then $(H^{\tau})_{o}=\mathrm{SO}_{o}(1,n-1)$.

\begin{lemma}
$\mathbb{R}^{+}\mathrm{SO}_{o}(1,n-1)$ is a symmetric subgroup of $H$.
\end{lemma}

We notice that $\frak{a}$ is maximal abelian in $\frak{s}$, $\frak
{s\cap q}$, and $\frak{q}$. Define $\alpha\in\frak{a}^{\ast}$ by 
$\alpha(tX)=t$ and $\alpha(\mathrm{I}_{n+1})=0$. Then $\Delta=\left\{ 
\alpha,-\alpha\right\} $. We choose $\alpha$ as the positive root. 
Then
\[
\frak{g}_{\alpha}=\left\{  X(v)=\left(
\begin{array}
[c]{ccc}
0 & v^{t} & 0\\
v & 0_{n-1} & -v\\
0 & v^{t} & 0
\end{array}
\right)  \mid v\in\mathbb{R}^{n-1}\right\}
\]
and
\[
N=\left\{  n(v)=\left(
\begin{array}
[c]{ccc}
1+\frac{1}{2}\left|  \left|  v\right|  \right|  ^{2} & v^{t} & -\frac{1}
{2}\left|  \left|  v\right|  \right|  ^{2}\\
v & \mathrm{I}_{n-1} & -v\\
\frac{1}{2}\left|  \left|  v\right|  \right|  ^{2} & v^{t} & 1-\frac{1}
{2}\left|  \left|  v\right|  \right|  ^{2}
\end{array}
\right)  \mid v\in\mathbb{R}^{n-1}\right\}  \,.
\]
In this case $\frak{m}_{2}=\left\{0\right\}$ and
\[
\frak{m}=\frak{b}\simeq\frak{so}(n-1)\ni X\mapsto\left(
\begin{array}
[c]{ccc}
0 & 0 & 0\\
0 & X & 0\\
0 & 0 & 0
\end{array}
\right)  \subset\frak{so}(n).
\]
The group $M$ is given by
\[
M\simeq\mathrm{SO}_{o}(n-1)\ni k\mapsto\left(
\begin{array}
[c]{ccc}
1 & 0 & 0\\
0 & k & 0\\
0 & 0 & 1
\end{array}
\right)  \subset\mathrm{SO}_{n}(n)\subset\mathrm{SO}_{o}(1,n)\,
\]
and
\[
A=\left\{  a(\lambda,t)=\lambda\left(
\begin{array}
[c]{ccc}
\cosh(t) & 0 & \sinh(t)\\
0 & \mathrm{I}_{n-1} & 0\\
\sinh(t) & 0 & \cosh(t)
\end{array}
\right)  \mid t\in\mathbb{R}\text{,}\lambda>0\right\}.
\]
Each of the groups $M$, $A$, and $N$ act on $\mathbb{R}^{n}$ in the 
following way:

\begin{itemize}
\item $M$: Rotation in the $v_{2},\ldots,v_{n-1}$ coordinates;

\item $A$: Dilations and hyperbolic rotations:
\[
a(\lambda,t)v=\lambda(\cosh(t)v_{1}+\sinh(t)v_{n+1},v_{2},\ldots,v_{n}
,\sinh(t)v_{1}+\cosh(t)v_{n+1})^{t}
\]

\item $N$: Write a vector $v\in\mathbb{R}^{n+1}$ as $v=(a,x,b)$ with
$a,b\in\mathbb{R}$ and $x\in\mathbb{R}^{n-1}$. Then
\[
n(u)v=v+(a-b)\left(
\begin{array}
[c]{c}
\frac{1}{2}\left|  \left|  u\right|  \right|  ^{2}\\
u\\
\frac{1}{2}\left|  \left|  u\right|  \right|  ^{2}
\end{array}
\right)  +(x\mid v)\left(
\begin{array}
[c]{c}
1\\
0\\
1
\end{array}
\right)
\]
\end{itemize}

The Weyl group $W$ consists of two elements $W=\left\{1,w\right\}$ 
where $w$ is given by $X\rightarrow-X$. It is realized by conjugating 
by $s$ where $s=\left(
\begin{array}
[c]{cccc}
1 &  &  & \\
& -1 &  & \\
&  & I_{n-2} & \\
&  &  & -1
\end{array}
\right)  $. Notice that $s\in H^{\tau}\setminus(H^{\tau})_{o}$ by
(\ref{eqtau}), and that $s(e_{n+1})=-e_{n+1}$. It follows that 
$\mathcal{O}_{3}$ decomposes into two open dense $P$-orbits 
$\mathcal{O}_{3,1}$, $\mathcal{O}_{3,2}$ and orbits of lower 
dimension.

\begin{theorem}
Let $n(v)\in N$, $a(\lambda,t)\in A$, and $k\in M$. Then
\begin{align*}
n(v)a(\lambda,t)k\cdot e_{n+1}=
\lambda\left(\sinh(t)-\frac{e^{-t}}{2}\left|\left|v\right|\right|^{2},
-e^{-t}v,\cosh(t)-\frac{e^{-t}}{2}\left|\left|v\right|\right|^{2}\right)^{t}\\
n(v)a(\lambda,t)k\cdot(-e_{n+1})
=\lambda\left(-\sinh(t)+\frac{e^{-t}}{2}\left|\left|v\right|\right|^{2},
e^{-t}v,-\cosh(t)+\frac{e^{-t}}{2}\left|\left|v\right|\right|^{2}\right)^{t}
\end{align*}
In particular $\mathcal{O}_{3,2}=-\mathcal{O}_{3,1}$ for
\begin{align*}
\mathcal{O}_{3,1}=Pe_{n+1}=\left\{  (a,v,b)\mid\left|  \left|  v\right|
\right|  ^{2}>a^{2}-b^{2}\,\,and\,\,a<b\right\}\text{\rm and } \\
\mathcal{O}_{3,2}=P(-e_{n+1})=\left\{  (a,v,b)\mid\left|  \left|  v\right|
\right|  ^{2}>a^{2}-b^{2}\,\,and\,\,a>b\right\}  .
\end{align*}
\end{theorem}

To avoid the problem of having to change the group as one goes from 
the orbits $\mathcal{O}_{1}$ and $\mathcal{O}_{2}$ to the orbit 
$\mathcal{O}_{3}$ we notice that $H=PK$, where $K=\mathrm{SO}_{o}(n)
=H^{\pm e_{1}}$. Hence $P$ actually acts transitively on the orbits 
$\mathcal{O}_{1}$ and $\mathcal{O} _{2}$ and in both cases we have a 
compact stabilizers. Thus

\begin{theorem}
Under the action of the group $P\times_{\pi}\mathbb{R}^{1+n}$, the 
space $L^{2}(\mathbb{R}^{1+n})$ decomposes into four irreducible parts
\[
L^{2}(\mathbb{R}^{1+n})\simeq L^{2}(\mathcal{O}_{1})\oplus L^{2}
(\mathcal{O}_{2})\oplus L^{2}(\mathcal{O}_{3,1})\oplus L^{2}(\mathcal{O}
_{3,2})\,.
\]
If we replace the central subgroup $\mathbb{R}^{+}$ by $\mathbb{R}^
{\ast}$, then there are only two irreducible parts:
\[
L^{2}(\mathbb{R}^{1+n})\simeq L^{2}(\mathcal{O}_{1}\cup\mathcal{O}_{2})\oplus
L^{2}(\mathcal{O}_{3,1}\cup\mathcal{O}_{3,2}).
\]
All the irreducible subrepresentations are square integrable, i.e., 
allow for dense subspaces of wavelet vectors.
\end{theorem}

\vskip 5pt
{\sc \bf\large Example} $
(\mathbb{R}^{+}\mathrm{SL(n,}\mathbb{R}),\mathrm{Symm}(n,\mathbb{R}))$:
\vskip 3pt

In the last examples in turned out, that one could
do the continuous wavelet transform uniformly by
replacing the group $H=\mathbb{R}^+\mathrm{SO}_o(n,\mathbb{R})$
by the parabolic subgroups $P$. We include the next example
to show that this is in general \textit{not} the case, and
in particular to show that the group $R$ used in the reconstruction
of $f$ from its wavelet transform, depends in general on
the orbit $\mathcal{O}_j$. 

Let $V$ be the $\frac{n(n+1)}{2}$-dimensional Euclidean vector space
\[
\mathrm{Symm}(n,\mathbb{R})=\left\{  X\in M_{n}(\mathbb{R})\mid X^{t}
=X\right\}
\]
with the inner product $(X,Y)=\mathrm{Tr}(XY^{t})$. The group 
$H=\mathbb{R} ^{+}\mathrm{SL}(n,\mathbb{R)}$ acts on $V$ by
\[
g\cdot X=gXg^{t}.
\]
The orbits are open if they contain nondegenerate 
bilinear forms and in this case are parameterized by their signature.
Specifically, let $p,q\in\mathbb{N}_{0}$ be 
such that $p+q=n$. Set $\mathrm{I}_{p,q}=\left(
\begin{array}
[c]{cc}
\mathrm{I}_{p} & 0\\
0 & -\mathrm{I}_{q}
\end{array}
\right)  $ and let
\[
\mathcal{O}_{p,q}=H\cdot\mathrm{I}_{p,q}=\left\{\text{symmetric 
matricies of signature }(p,q)\right\}\,.
\]
Take $\mathrm{SO}(p,q)=\left\{  g\in\mathrm{SL}(n,\mathbb{R})\mid
g\mathrm{I}_{p,q}g^{t}=\mathrm{I}_{p,q}\right\}$, and define $\tau
_{p,q}:H\rightarrow H$ by 
\[\tau_{p,q}(\lambda g)=
\lambda^{-1}\mathrm{I}_{p,q}(g^{t})^{-1}\mathrm{I}_{p,q}, 
\quad\lambda>0,\quad g\in\mathrm{SL}(n,\mathbb{R}).
\] 
Then $\tau_{p,q}:H\rightarrow H$ is an involution and
$H^{\tau}=\mathrm{SO}(p,q)$. The following is now clear:

\begin{theorem}
Let the notation be as above. Then the following hold:

\begin{enumerate}
\item \vskip -.2cm 
Each orbit $\mathcal{O}_{p,q}$ is open in $V$.

\item 
$\mathcal{O}_{p,q}=H/\mathrm{SO}(p,q)$. In particular it follows that 
the stabilizer of a point in $\mathcal{O}_{p,q}$ is compact if and 
only if $n=p$ or $n=q$.

\item 
$V\setminus\bigcup_{p+q=n}\mathcal{O}_{p,q}=\left\{ X\in V\mid
\det(X)=0\right\}  $. In particular it follows that $\bigcup_{p+q=n}
\mathcal{O}_{p,q}$ is dense and open in $V$.
\end{enumerate}
\end{theorem}

We note that the parabolics $P_{p,q}$ and the corresponding groups 
$R_{p,q}$ depend on the orbits $\mathcal{O}_{p,q}$.
These can be worked out in detail but we leave the details to
the reader.

\end{document}